\magnification=1200
\overfullrule=0pt
\centerline {{\bf Sublevel sets and global minima of coercive functionals}}
\centerline  {{\bf and local minima
of their perturbations}}\par
\bigskip
\bigskip
\centerline {BIAGIO RICCERI}\par
\bigskip
\bigskip
\bigskip
\bigskip
\noindent
{\bf Introduction.}\par
\bigskip
In [24],  we have identified a general variational principle of which the following
theorem is a by-product:\par
\medskip
THEOREM A. - {\it
 If $\Phi$ and $\Psi$ are two sequentially weakly lower semicontinuous functionals on a
reflexive real Banach space and if $\Psi$ is also continuous and coercive,
then the functional $\Psi+\lambda \Phi$ has at least a local minimum for each
$\lambda>0$ small enough.}\par
\medskip
The variational principle of [24] has already been widely applied
to nonlinear differential equations  and
Hammerstein integral equations as well (see, for instance, [1]-[6],
[8]-[16], [18], [20], [25]-[28]).\par
\smallskip
The aim of the present paper is essentially to point out that,
under the same assumptions as those of Theorem A, the following
 more precise conclusion holds:\par
\smallskip
 if, for some
$r>\inf_{X}\Psi$, the weak closure of the set $\Psi^{-1}(]-\infty,r[)$
has at least $k$ connected components in the weak topology, then, for
each $\lambda>0$ small enough, the functional $\Psi+\lambda\Phi$ has
at least $k$ local minima lying in $\Psi^{-1}(]-\infty,r[)$.\par
\smallskip
 This,
in particular, holds (for any $r>\inf_{X}\Psi$) when the set of all global
 minima of $\Psi$ has at least $k$ connected components in the
weak topology.\par
\smallskip
This more precise conclusion can be used in a twofold way.\par
\smallskip
 In a direct way, 
when we have, {\it a priori},  a sufficient information 
about the set of all global minima or, more generally, about the sublevel sets
of $\Psi$, it just provides an information on the number of  the
local minima of suitable
perturbations of $\Psi$. \par
\smallskip
Otherwise, when our primary objective is to get some
information on the structure of the set of all global minima and of the
sublevel
sets of $\Psi$, we can try to use it in an indirect way.\par
\smallskip
 For instance, if  we
are interested in knowing whether the sublevel sets of $\Psi$
 are connected in the weak topology (an important issue in
 minimax theory (see [22], [23])), then
we could try to find 
 a sequentially weakly lower semicontinuous functional $\Phi$ and
a sequence of positive numbers $\{\mu_{n}\}$ converging to $0$ in
such a way that, for each $n\in {\bf N}$, the functional $\Psi+\mu_{n}\Phi$
has at most one local minimum.\par
\smallskip
We develop this point of view in the third section, when $\Psi$ is the energy
functional related to a Dirichlet problem. \par
\smallskip
In the next section, we first establish our basic results in full generality and then
we formulate them in the setting of reflexive and separable real Banach spaces.
\vfill\eject
\noindent
{\bf Basic results.}\par
\bigskip
If $(X,\tau)$ is a topological space, for any $\Psi:X\to ]-\infty,+\infty]$,
with dom$(\Psi)\neq \emptyset$, we denote by $\tau_{\Psi}$ the smallest
topology on $X$ which contains both $\tau$ and the family of sets
$\{\Psi^{-1}(]-\infty,r[)\}_{r\in {\bf R}}$.\par
\smallskip
Our main abstract result is as follows.\par
\medskip
THEOREM 1. - {\it Let $(X,\tau)$ be a
Hausdorff topological space, and $\Phi, \Psi:X\to
]-\infty,+\infty]$ two functions.
Assume that there is $\rho>\inf_{X}\Psi$ such that the set
 $\overline {\Psi^{-1}(]-\infty,\rho[)}$ is compact and sequentially
 compact,
 has at least $k$ connected components and 
each of them  intersects the interior of
 \hbox {\rm dom}$(\Phi)$.
Moreover, suppose that the function $\Phi$ is
bounded below in $\overline {\Psi^{-1}(]-\infty,\rho[)}$ and that the
function $\Psi+\lambda\Phi$ is sequentially
lower semicontinuous for each $\lambda>0$ small enough.\par
Then, there exists $\lambda^{*}>0$ such that, for each
$\lambda\in ]0,\lambda^{*}[$, the function
$\Psi+\lambda\Phi$ has at least $k$ $\tau_{\Psi}$-local minima lying
in \hbox {\rm dom}$(\Phi)\cap \Psi^{-1}(]-\infty,\rho[)$.}
\smallskip
PROOF.
Denote by ${\cal C}$ the family of all connected components of
$\overline {\Psi^{-1}(]-\infty,\rho[)}$.
 Note that these sets are closed in
$X$ since they are closed in
 $\overline {\Psi^{-1}(]-\infty,\rho[)}$ which is, in turn, closed in $X$.
We now observe that there
are $k$ pairwise disjoint closed non-empty sets
$C_{1},...,C_{k}$ 
 such that
$$\overline {\Psi^{-1}(]-\infty,\rho[)}=\bigcup_{i=1}^{k} C_i\ .$$
We distinguish two cases. First, assume that ${\cal C}$ is finite.
Let $h$ be its cardinality. Let $B_{1},..,B_{h}$ be
 the members of ${\cal C}$. Then, if we choose $C_{i}=B_{i}$ for
$i=1,...,k-1$ and $C_{k}=\cup_{i=k}^{h}B_{i}$, we are clearly done.
Now, assume that ${\cal C}$ is infinite. In this case, we prove our claim
by induction. The claim is true, of course, if $k=1$.
Assume that it is true if $k=p$. So, we are assuming
that there are $p$ pairwise disjoint closed non-empty sets $D_{1},...,D_{p}$,
such that
$$\overline {\Psi^{-1}(]-\infty,\rho[)}=\bigcup_{i=1}^{p} D_i\ .$$
Notice that at least one of the sets $D_i$ must be disconnected, since,
otherwise, we would have
$\{D_{1},...,D_{p}\}={\cal C}$, contrary to the assumption that ${\cal C}$
is infinite. 
Then, if $D_{i^*}$ is disconnected, there are two disjoint closed non-empty
sets $E_{1}, E_{2}$
such that $D_{i^*}=E_{1}\cup E_{2}$. So,
$D_{1},...,D_{i^{*}-1},D_{i^{*}+1},...,D_{p},E_{1},E_{2}$ are
$p+1$ pairwise disjoint closed non-empty sets whose union
is $\overline {\Psi^{-1}(]-\infty,\rho[)}$. So, our claim
is true for $k=p+1$, and hence, by induction, for any $k$.\par
Now, fix $i$ ($1\leq i\leq k$). By compactness and Hausdorffness,
 it is clear that there exists
an open set $A_{i}\subset X$ such that $C_{i}\subset A_{i}$ and
$A_{i}\cap \cup_{j=1,j\neq i}^{k}C_{j}=\emptyset$. Furthermore, it is easily
seen that, if we put
$$G_{i}=\{x\in A_{i} : \Psi(x)<\rho\}\ ,$$
we have
$$\overline {G_{i}}=C_{i}\ .$$
 Clearly, $C_{i}$ contains at least one
member of ${\cal C}$, and hence, by assumption,
int(dom$(\Phi)$)$\cap C_{i}\neq \emptyset$. This implies that
dom$(\Phi)\cap G_{i}\neq \emptyset$.
Taken into account that, by assumption,
$\inf_{C_{i}}\Phi$ is finite, put
$$\mu_{i}=\inf_{x\in \hbox {\rm dom}(\Phi)\cap G_{i}}
{{\Phi(x)-\inf_{C_{i}}\Phi}\over {\rho-\Psi(x)}}\ .$$
Let $\lambda '>0$ be such that
$\Psi+\lambda\Phi$ is sequentially lower semicontinuous for each
$\lambda\in ]0,\lambda ']$.
Fix $\mu>\max\{\mu_{i},{{1}\over {\lambda '}}\}$. Then,
there exists $y\in \hbox {\rm dom}(\Phi)
\cap G_{i}$ such that
$$\mu\rho>\mu\Psi(y)+\Phi(y)-\inf_{C_{i}}\Phi\ .$$
Moreover, since $C_i$ is sequentially compact,
 there exists $x^{*}_{i}\in \hbox {\rm dom}(\Phi)
\cap C_{i}$ such
$$\Phi(x^{*}_{i})+\mu\Psi(x^{*}_{i})\leq \Phi(x)+\mu\Psi(x)$$
for all $x\in C_{i}$. We claim that $x^{*}_{i}\in G_{i}$. Arguing by
contradiction, assume that $\Psi(x^{*}_{i})\geq \rho$. We then have
$$\Phi(x^{*}_{i})+\mu\Psi(x^{*}_{i})\geq
\Phi(x^{*}_{i})+\mu\rho>\Phi(x^{*}_{i})+\Phi(y)+\mu\Psi(y)-\inf_{C_i}\Phi\geq
\Phi(y)+\mu\Psi(y)$$
which is absurd. Now, let $i$ vary. Put $\mu^{*}=
\max\{\mu_{1},...,\mu_{k},{{1}\over {\lambda '}}\}$.
Clearly, each set $G_{i}$ is $\tau_{\Psi}$-open, and 
hence each $x^{*}_{i}$ is a $\tau_{\Psi}$-local minimum of $\Phi+\mu\Psi$ for all
$\mu>\mu^{*}$. Consequently, the points $x^{*}_{1},...,x^{*}_{k}$ satisfy the conclusion,
taking $\lambda^{*}={{1}\over {\mu^{*}}}$, and the proof is complete.\hfill 
$\bigtriangleup$\par
\medskip
The next result provides a reasonable way to check the key assumption of
Theorem 1.\par
\medskip
PROPOSITION 1. - {\it Let $X$ be a Hausdorff topological space and
$\Psi:X\to ]-\infty,+\infty]$ a sequentially
lower semicontinuous function. Assume that there is $r>\inf_{X}\Psi$
such that the set
$\overline {\Psi^{-1}(]-\infty,r[)}$ is compact and first-countable. Moreover,
assume that
the set of all global minima of $\Psi$ has at least $k$ connected
components.\par
Then, there exists $\rho^{*}\in ]\inf_{X}\Psi,r]$ such that, for each
$\rho\in ]\inf_{X}\Psi,\rho^{*}]$,
the set $\overline {\Psi^{-1}(]-\infty,\rho[)}$ has at least $k$ connected
components.}\par
\smallskip
PROOF. Arguing by contradiction, assume that there is a decreasing sequence
$\{\rho_{n}\}$ in $]-\infty,r[$, coverging to $\inf_{X}\Psi$, such that
for each $n\in {\bf N}$, the set $\overline {\Psi^{-1}(-\infty,\rho_{n}[)}$
has at most $k-1$ connected components. Clearly, we have
$$\Psi^{-1}(\inf_{X}\Psi)=\bigcap_{n=1}^{\infty}
\overline {\Psi^{-1}(]-\infty,\rho_{n}[)}\ .$$
Reasoning as in the proof of
Thoerem 1, we find $k$ open and pairwise disjoint subsets of $X$,
 $\Omega_{1},...,\Omega_{k}$, such that
 $\Omega_{i}\cap \Psi^{-1}(\inf_{X}\Psi)\neq \emptyset$ for all
 $i=1,...,k$ and $\Psi^{-1}(\inf_{X}\Psi)\subseteq \bigcup_{i=1}^{k}\Omega_i$.
Then, for each $n\in {\bf N}$, the set
$\overline {\Psi^{-1}(]-\infty,\rho_{n}[)}$ cannot be contained in
$\bigcup_{i=1}^{k}\Omega_i$, since, otherwise, it would have at
least $k$ connected components. So, 
$\overline {\Psi^{-1}(]-\infty,\rho_{n}[)}\cap 
(X\setminus \bigcup_{i=1}^{k}\Omega_{i})$ is a non-increasing sequence
of non-empty closed subsets of a compact one, and hence 
$\bigcap_{n=1}^{\infty}
\overline {\Psi^{-1}(]-\infty,\rho_{n}[)}\cap 
(X\setminus \bigcup_{i=1}^{k}\Omega_{i})\neq \emptyset$
which is absurd.\hfill $\bigtriangleup$\par
\medskip
So, putting Theorem 1 and Proposition 1 together, we clearly get\par
\medskip
THEOREM 2. - {\it  Let $(X,\tau)$ be a
Hausdorff topological space and $\Phi:X\to {\bf R}$, 
$\Psi:X\to ]-\infty,+\infty]$ two
functions.
Assume that there is $r>\inf_{X}\Psi$ such that the set
$\overline {\Psi^{-1}(]-\infty,r[)}$ is compact and first-countable.
Moreover, suppose that the function $\Phi$ is
bounded below in $\overline {\Psi^{-1}(]-\infty,r[)}$ and that the
function $\Psi+\lambda\Phi$ is sequentially
lower semicontinuous for each $\lambda\geq 0$ small enough.
 Finally,
assume that
the set of all global minima of $\Psi$ has at least $k$ connected
components.\par
Then, there exists $\lambda^{*}>0$
such that, for each
$\lambda\in ]0,\lambda^{*}[$, the function $\Psi+\lambda\Phi$ has
at least $k$ $\tau_{\Psi}$-local minima
lying in $ \Psi^{-1}(]-\infty,r[)$.}\par
\medskip
Arguing by contradiction, the use of Theorems 1 and 2 gives the following
\medskip
THEOREM 3. -
{\it Let $(X,\tau)$ be a Hausdorff topological space and
$\Psi:X\to ]-\infty,+\infty]$ a sequentially lower semicontinuous function
such that, for some $r>\inf_{X}\Psi$, the set
$\overline {\Psi^{-1}(]-\infty, r[)}$ is compact and first-countable. \par
  Suppose that there
are a function $\Phi:X\to {\bf R}$, bounded below in
$\overline {\Psi^{-1}(]-\infty, r[)}$,
  and 
 a sequence $\{\mu_{n}\}$ in ${\bf R}^+$ 
converging to $0$ such that, for each $\lambda>0$ small enough, the function
$\Psi+\lambda\Phi$ is sequentially lower semicontinuous, and, for
each $n\in {\bf N}$, the function $\Psi+\mu_{n}\Phi$ has at most $k$
$\tau_{\Psi}$-local minima lying in $\Psi^{-1}(]-\infty,r[)$.\par
Then, for every $\rho\in ]\inf_{X}\Psi, r]$,  the sets
$\overline {\Psi^{-1}(]-\infty, \rho[)}$ and
$\Psi^{-1}(\inf_{X}\Psi)$ have at most $k$ connected components.
So, in particular, these sets are connected when $k=1$.}\par
\medskip
REMARK 1. - When $k=1$, Theorem 3 ensures that, for every
$\rho\in ]\inf_{X}\Psi, r[$, the set 
$\Psi^{-1}(]-\infty, \rho])$ is connected.
This follows from the equality
$$\Psi^{-1}(]-\infty,\rho])=\bigcap_{\rho<s<r}
\overline {\Psi^{-1}(]-\infty,s[)}$$
and from the fact that, for every $s\in ]\inf_{X}\Psi,r]$, the
set  $\overline {\Psi^{-1}(]-\infty,s[)}$ is connected and compact.
\medskip
An interesting consequence of Theorem 3 is the following
two local minima result.\par
\medskip
THEOREM 4. - {\it  Let $(X,\tau)$ be a
Hausdorff topological space, and $\Phi:X\to {\bf R}$, 
$\Psi:X\to ]-\infty,+\infty]$ two functions.
Assume that there is $r>\inf_{X}\Psi$ such that the set
$\overline {\Psi^{-1}(]-\infty,r[)}$ is compact and first-countable.
 Moreover,
assume that there is a strict local minimum of $\Psi$, say $x_0$,  such that
$\inf_{X}\Psi<\Psi(x_{0})<r$.
Finally, suppose that the function $\Phi$ is
bounded below in $\overline {\Psi^{-1}(]-\infty,r[)}$ and that the
function $\Psi+\lambda\Phi$ is sequentially
lower semicontinuous for each $\lambda\geq 0$ small enough.\par
Then, there exists $\lambda^{*}>0$ such that, for each
$\lambda\in ]0,\lambda^{*}[$, the function $\Psi+\lambda\Phi$ has
at least two $\tau_{\Psi}$-local minima
lying in $ \Psi^{-1}(]-\infty,r[)$.}\par
\smallskip
PROOF. Arguing by contradiction, assume that the conclusion does not hold.
Then, by  Theorem 3 and Remark 1,
the set $\Psi^{-1}(]-\infty,\Psi(x_{0})])$ is connected. But this set
contains $x_{0}$ as an isolated point (since $x_0$ is a strict local
minimum of $\Psi$) and does not reduce to it (since $\inf_{X}\Psi<
\Psi(x_{0})$), against connectedness.\hfill $\bigtriangleup$\par
\medskip
With the aim to apply them to nonlinear differential equations,
we now establish some consequences of the previous general
results in the setting of reflexive and separable real Banach spaces. \par
\smallskip
For a set $A$ in a Banach space, we denote by $(\overline {A})_w$ its
closure in the weak topology. We say that $A$ is weakly connected if it
is connected in the weak topology. The weakly connected components of $A$
are its connected components in the weak topology.
\medskip
THEOREM 5. -
{\it Let $X$ be a sequentially weakly closed
subset of a reflexive and separable real Banach space $E$,
 and  $\Phi:X\to {\bf R}$,
$\Psi:X\to ]-\infty,+\infty]$ two  
functionals.
Assume that there is
$\rho>\inf_{X}\Psi$ such that the set
$(\overline {\Psi^{-1}(]-\infty, \rho[)})_w$ is bounded and has
at least $k$ weakly connected components. Moreover,
suppose that the functional $\Phi$ is bounded
below in $(\overline {\Psi^{-1}(]-\infty, \rho[)})_w$
and that the functional
$\Psi+\lambda\Phi$ is
sequentially weakly lower semicontinuous for each $\lambda>0$
small enough.\par
Then, there exists $\lambda^{*}>0$ such that, for each
$\lambda\in ]0,\lambda^{*}[$, the functional
$\Psi+\lambda\Phi$ has at least $k$ $\tau_\Psi$-local minima lying
in  $\Psi^{-1}(]-\infty,\rho[)$, where $\tau$ is the relative weak topology
of $X$.}\par
\smallskip
PROOF. Apply Theorem 1, $\tau$ just being the relative weak topology.
In particular,
observe that, since $E$ is reflexive and separable, the weak closure of
any bounded set is weakly compact and metrizable (and so first-countable). 
\hfill $\bigtriangleup$\par
\medskip
Analogously, from Theorem 2, we get\par
\medskip
THEOREM 6. -
{\it Let $X$ be a sequentially weakly closed
subset of a reflexive and separable real Banach space $E$,
 and  $\Phi:X\to {\bf R}$,
$\Psi:X\to ]-\infty,+\infty]$ two  
functionals.
Assume that there is
$\rho>\inf_{X}\Psi$ such that the set
$\Psi^{-1}(]-\infty, \rho[)$ is bounded. Moreover,
suppose that the functional $\Phi$ is bounded         
below in $(\overline {\Psi^{-1}(]-\infty,\rho[)})_w$
and that the functional
$\Psi+\lambda\Phi$ is
sequentially weakly lower semicontinuous for each $\lambda\geq 0$
small enough. Finally,
assume that
the set $\Psi^{-1}(\inf_{X}\Psi)$ has at least $k$ weakly connected
components.\par
Then, the conclusion of Theorem 5 holds}.\par
\medskip
Arguing by contradiction, from Theorems 5 and 6 we then get
\medskip
THEOREM 7. - 
{\it Let $X$ be a sequentially weakly closed
subset of a reflexive and separable real Banach space, and
$\Psi:X\to ]-\infty,+\infty]$ a sequentially weakly
lower semicontinuous functional such that, for some
$r>\inf_{X}\Psi$, the set
$\Psi^{-1}(]-\infty, r[)$ is bounded.\par
  Suppose that there
are a functional $\Phi:X\to {\bf R}$, bounded below in
$(\overline {\Psi^{-1}(]-\infty, r[)})_{w}$,
  and 
 a sequence $\{\mu_{n}\}$ in ${\bf R}^+$ 
converging to $0$ such that, for each $\lambda>0$ small enough, the
functional
$\Psi+\lambda\Phi$ is sequentially weakly lower semicontinuous, and, for
each $n\in {\bf N}$, the functional $\Psi+\mu_{n}\Phi$ has at most $k$
$\tau_{\Psi}$-local minima lying in $\Psi^{-1}(]-\infty,r[)$,
where $\tau$ is relative weak topology of $X$\par
Then, for every $\rho\in ]\inf_{X}\Psi, r]$,  the sets
$\overline {\Psi^{-1}(]-\infty, \rho[)}$ and
$\Psi^{-1}(\inf_{X}\Psi)$ have at most $k$ weakly connected components.
So, in particular, these sets are weakly connected when $k=1$.}\par
\medskip
The next result is an application of Theorems 5 and 6 to critical point
theory. If $J$ is a G\^ateaux differentiable functional on a Banach
space $X$, the critical points of $J$ are the zeros of its derivative, $J'$.
Moreover, $J$ is said to satisfy the Palais-Smale condition  if
each sequence $\{x_{n}\}$ in $X$ such that $\sup_{n\in {\bf N}}|
J(x_{n})|<+\infty$ and $\lim_{n\to +\infty}\|J'(x_{n})\|_{X^{*}}=0$ admits
a strongly converging subsequence.\par
\medskip
THEOREM 8. - {\it In addition to the assumptions of either
Theorem 5 or Theorem 6, suppose
that $X=E$, that the functionals
$\Psi, \Phi: X\to {\bf R}$ are
continuously G\^ateaux differentiable, and that $k\geq 2$.\par 
Then, there exists $\lambda^{*}>0$
such that, for each
$\lambda\in ]0,\lambda^{*}[$ for which the functional $\Psi+\lambda\Phi$
satisfies
the Palais-Smale condition, the same functional has
at least $k+1$ critical points, $k$ of which are lying in
$ \Psi^{-1}(]-\infty,\rho[)$.}\par
\smallskip
PROOF. By either Theorem 5 or Theorem 6,  there exists $\lambda^{*}>0$ such that, for each
$\lambda\in ]0,\lambda^{*}[$, the functional
$\Psi+\lambda\Phi$ has at least $k$ $\tau_\Psi$-local minima lying
in  $\Psi^{-1}(]-\infty,\rho[)$, where $\tau$ is the weak topology
of $X$. Note that $\Psi$, being $C^1$, is (norm) continuous. Consequently,
the topology $\tau_\Psi$ is weaker than the strong topology, and so
the above mentioned $\tau_\Psi$-local minima of $\Psi+\lambda\Phi$ are
local minima of this functional in the strong topology. Now, assuming
that $\Psi+\lambda\Phi$ satisfies the Palais-Smale condition,
the conclusion follows from Theorem (1.ter) of [17].\hfill $\bigtriangleup$\par
\medskip
 From Theorem 8, arguing by contradiction, we then obtain the following\par
\medskip
THEOREM 9. - {\it Let $X$ be a reflexive and separable real Banach space
and $\Psi : X\to {\bf R}$ a continuously G\^ateax differentiable
 and sequentially weakly lower semicontinuous
functional such that, for some $r>\inf_{X}\Psi$, the set
$\Psi^{-1}(]-\infty,r[)$ is bounded. Let
$k\in {\bf N}$ with $k\geq 2$.\par
  Suppose that there are a continuously
G\^ateaux differentiable functional
 $\Phi : X\to {\bf R}$, which is bounded below in
$(\overline {\Psi^{-1}(]-\infty, r[)})_w$, and
a sequence $\{\mu_{n}\}$ in ${\bf R}^+$ 
converging to $0$ such that the functional
$\Psi+\lambda\Phi$ is sequentially weakly lower semicontinuous
for each $\lambda>0$ small enough and, for
each $n\in {\bf N}$, the functional $\Psi+\mu_{n}\Phi$
satisfies the Palais-Smale condition and has at most $k$
critical points in $X$.\par
Then, for every $\rho\in ]\inf_{X}\Psi, r]$, the sets
$(\overline {\Psi^{-1}(]-\infty, \rho[)})_w$ and
$\Psi^{-1}(\inf_{X}\Psi)$ have at most $k-1$ weakly connected components.
So, in particular, these sets are weakly connected when $k=2$.}\par
\medskip 
REMARK 2. - When $k=2$, Theorem 9 ensures that, for every
$\rho\in ]\inf_{X}\Psi, r[$, the set 
$\Psi^{-1}(]-\infty, \rho])$ is weakly connected (see Remark 1).\par
\medskip
Here is an application of Theorem 4.\par
\medskip
THEOREM 10. - {\it Let $X$ be a uniformly convex and separable real Banach
space, $g : [0,+\infty[\to {\bf R}$ a strictly increasing continuous function,
and $J : X\to {\bf R}$ a sequentially weakly lower
semicontinuous functional. For every $x\in X$, put
$$\Psi(x) = g(\|x\|) + J(x)\ .$$
Assume that the functional $\Psi$ is coercive and has a strict,
not global, local minimum, say $x_0$.\par
Then,  for every 
$r>\Psi(x_{0})$ and every functional $\Phi:X\to {\bf R}$ which is
bounded in $(\overline {\Psi^{-1}(]-\infty, r[)})_w$ and such that
$\Psi+\lambda\Psi$ is sequentially weakly lower semicontinuous for each
$\lambda>0$ small enough,
there exists $\lambda^{*}>0$ such that, for each
$\lambda\in ]0,\lambda^{*}[$, the functional  $\Psi + \lambda\Phi$ has
at least two $\tau_{\Psi}$-local minima lying in $ \Psi^{-1}(]-\infty,r[)$,
where $\tau$ is the weak topology of $X$}.\par
\smallskip
PROOF. From Theorem 1 of [21], it follows that $x_{0}$ is a $\tau$-strict local
minimum of $\Psi$. More precisely, the statement of the above quoted result
deals with local minima, but exactly the same proof shows that the same is
true for strict local minima. Now, the conclusion follows
from Theorem 4, taking as $\tau$ just the weak topology.\hfill
$\bigtriangleup$\par
\medskip
We conclude this section with an application of Theorem 7 in the setting
of Hilbert spaces.\par
\medskip
THEOREM 11. - {\it 
Let $X$ be a separable real Hilbert space and
$J: X\to {\bf R}$ a continuous,
 G\^ateaux differentiable 
 and sequentially weakly upper
semicontinuous functional.
For every $x\in X$, put
$$\Psi(x) = {{1}\over {2}}\|x\|^{2} - J(x)\ .$$
Assume that, for some $r>\inf_{X}\Psi$, the set $\Psi^{-1}(]-\infty,r[)$ is
bounded. Moreover, suppose that the restriction of $J'$ to 
$\Psi^{-1}(]-\infty,r[)$ is nonexpansive.\par
Then, for every $\rho\in ]\inf_{X}\Psi, r[$, the sets
$ \Psi^{-1}(]-\infty, \rho])$ and
$\Psi^{-1}(\inf_{X}\Psi)$ are weakly connected. }
\smallskip
PROOF. Let us apply Theorem 7 taking 
as $\{\mu_{n}\}$ any sequence in $]0,+\infty[$ coverging
to zero, and $\Phi(x)={{1}\over {2}}\|x\|^2$.
Thus, $\Phi$ and $\Psi$ are
two continuous, G\^ateaux differentiable and
sequentially weakly lower semicontinuous functionals, and,
 for each $n\in {\bf N}$, the functional $\Psi+\mu_{n}\Phi$ 
 admits
 at most one $\tau_{\Psi}$-local minimum in $\Psi^{-1}(]-\infty,r[)$.
This follows from the fact that the strong topology is stronger than
$\tau_{\Psi}$ and that the
retriction of ${{1}\over {1+\mu_{n}}}J'$ to
$\Psi^{-1}(]-\infty,r[)$
is a contraction, and so the equation $\Psi'(x)+\mu_{n}\Phi'(x)=0$
has at most one solution in $\Psi^{-1}(]-\infty,r[)$. Hence, the
hypotheses of Theorem 7 are satisfied, and the
conclusion follows from it and from Remark 1.\hfill $\bigtriangleup$\par
\medskip
The following proposition is an useful complement to both
Theorems 9 and 11.\par
\medskip
PROPOSITION 2. - {\it 
Let $X$ be a real Hilbert space and
$J: X\to {\bf R}$ a continuously G\^ateaux differentiable
functional whose derivative is compact.
For every $x\in X$, put
$$\Psi(x) = {{1}\over {2}}\|x\|^{2} - J(x)\ .$$
Assume that, for some $r>\inf_{X}\Psi$, the set $\Psi^{-1}(]-\infty,r[)$ is
bounded.\par
Then, the set $\Psi^{-1}(\inf_{X}\Psi)$ is compact.}\par
\smallskip
PROOF. The functional $\Psi$ is G\^ateaux differentiable and its critical
points are exactly the fixed points of $J'$.
 Let $B$ be a closed ball in $X$ containing $\Psi^{-1}(]-\infty,r[)$. Let
$\{x_{n}\}$ be a sequence of fixed points of $J'$ lying in $B$. Since
this sequence is bounded and $J'$ is compact, there is a subsequence
$\{x_{n_{k}}\}$ such that $J'(x_{n_{k}})$ converges to some $z\in B$.
Clearly, by continuity, $z$ is a fixed point of $J'$. So, the set
$\{x\in B : J'(x)=x\}$ is compact. Of course, it contains 
$\Psi^{-1}(\inf_{X}\Psi)$ which is closed since $\Psi$ is continuous, and
the conclusion follows.\hfill $\bigtriangleup$\par
\medskip
REMARK 3. - Note that when we can apply
Theorem 7 with $k=1$ to a functional $\Psi$ as in Proposition 2, then the set 
$\Psi^{-1}(\inf_{X}\Psi)$ is connected. This follows from the fact
in any compact subset of a Banach space the relative strong and weak
topologies coincide. This remark, in particular, applies to Theorem 9 (when
$k=2$) and  to Theorem 11.
\par
\vfill\eject
\noindent
{\bf Applications.}
\bigskip
In this section, we intend to present some applications of Theorems 9 and
11 to the energy functional related to the Dirichlet problem for a semilinear elliptic
equation.
\par
\smallskip
In the sequel, $\Omega$ will denote an open connected subset of
${\bf R}^h$ with sufficiently smooth boundary.\par
\medskip
Put $X=W^{1,2}_{0}(\Omega)$, and consider it with the usual norm
$\|u\|=(\int_{\Omega}|\nabla u(x)|^{2}dx)^{1\over 2}$.
If $h\geq 2$, we denote by ${\cal A}$ the class of all
Carath\'eodory functions $f:\Omega\times {\bf R}\to {\bf R}$ such that
$$\sup_{(x,\xi)\in \Omega\times {\bf R}}{{|f(x,\xi)|}\over
{1+|\xi|^q}}<+\infty\ ,$$
where  $0<q< {{h+2}\over {h-2}}$ if $h>2$ and $0<q<+\infty$ if
$h=2$. While, when $h=1$, we denote by ${\cal A}$  the class
of all Carath\'eodory functions $f:\Omega\times {\bf R}\to {\bf R}$ such
that, for each $r>0$, the function $x\to \sup_{|\xi|\leq r}|f(x,\xi)|$ belongs
to $L^{1}(\Omega)$.\par
\smallskip
For each  $f\in {\cal A}$ and $u\in X$, we put
$$J_{f}(u)=
\int_{\Omega}\left ( \int_{0}^{u(x)}f(x,\xi)d\xi\right ) dx$$
and
$$\Psi_{f}(u)={{1}\over {2}}\int_{\Omega}|\nabla u(x)|^{2}dx-J_{f}(u)\ .$$
\smallskip
So,  by classical results, 
the functional $J_{f}$ is (well defined and) continuously G\^ateaux differentiable on
$X$, its derivative is compact, and one has
$$\Psi_{f}'(u)(v)=\int_{\Omega}\nabla u(x)\nabla v(x)dx-
\int_{\Omega}f(x,u(x))v(x)dx$$
for all $u, v\in X$. 
 Hence, the critical points of $\Psi_{f}$ in $X$ are exactly
the weak solutions of  the Dirichlet problem
$$\cases {-\Delta u=
f(x,u)
 & in
$\Omega$\cr & \cr u_{|\partial \Omega}=0\ . \cr} $$
Recall also that if $\Psi_{f}$ is coercive, then it satisfies the Palais-Smale
condition (see, for instance, Example 38.25 of [29]).
\medskip
We denote by $\lambda_{1}$ the first eigenvalue of the problem
$$\cases {-\Delta u=\lambda u & in
$\Omega$\cr & \cr u_{|\partial \Omega}=0 \ .\cr}$$
Recall that
$\|u\|_{L^{2}(\Omega)}\leq \lambda_{1}^{-{{1}\over {2}}}\|u\|$
for all $u\in X$. So, if
$$\limsup_{|\xi|\to +\infty}{{\sup_{x\in \Omega}\int_{0}^{\xi}f(x,t)dt}
\over {\xi^2}}<{{\lambda_{1}}\over {2}}$$
then the functional $\Psi_{f}$ is coercive in $X$.\par
\smallskip
Let us state the following\par
\medskip
THEOREM 12. - {\it
Let $f : \Omega\times {\bf R}\to {\bf R}$ be a Carath\'eodory function belonging
to ${\cal A}$
 such that
$$\sup_{(\xi,\eta)\in {\bf R}^{2},\hskip 3pt
\xi\neq \eta}{{\sup_{x\in \Omega}|f(x,\xi)-f(x,\eta)|}\over  {|\xi-\eta|}}\leq
\lambda_{1}$$
and
$$\limsup_{|\xi|\to +\infty}{{\sup_{x\in \Omega}\int_{0}^{\xi}f(x,t)dt}\over {\xi^2}}
<{{\lambda_{1}}\over {2}}\ .$$
Then, the sublevel sets of $\Psi_{f}$ are weakly
connected, and the set of all global minima of $\Psi_{f}$ is
compact and connected.}
\smallskip
PROOF.  
Fix $u, v, w\in X$, with $\|w\|=1$.
 We have
$$|J_{f}'(u)(w)-J_{f}'(v)(w)|\leq 
\int_{\Omega}|f(x,u(x))-f(x,v(x))||w(x)|dx\leq
\lambda_{1}\|u-v\|_{L^{2}(\Omega)}\|w\|_{L^{2}(\Omega)}\leq
\|u-v\|\ ,$$
and hence
$$\|J_{f}'(u)-J_{f}'(v)\|\leq \|u-v\|\ ,$$
 that is $J'$ is nonexpansive in $X$.
 Moreover, $\Psi_{f}$ is coercive in $X$. Thus, the functionals
$J_{f}$ and $\Psi_{f}$ satisfy all the assumptions of Theorem 11, and the
conclusion follows from it, taking also into account Proposition 2 and
Remark 3.\hfill $\bigtriangleup$\par
\medskip
Let us segnalize an open question related to Theorem 12.\par
\medskip
PROBLEM 1. - Is there some $f$ satisfying the assumptions of Theorem 11
for which the
set of all global minima of the functional $\Psi_{f}$ is neither a
singleton nor a segment ? 
In particular, what happens when $f(\xi)=\lambda_{1}(\sin\xi + a) $, with $a>0$ ? Or
when $f(\xi)=\lambda_{1}\hbox {\rm dist}(\xi,A)$, where $A\subset {\bf R}$ ?
\medskip
We now establish\par
\medskip
THEOREM 13. - {\it Let $g:\overline {\Omega}\times [0,+\infty[\to {\bf R}$
be a locally H\"older continuous function belonging to ${\cal A}$ such that
$$\limsup_{\xi\to +\infty}{{\sup_{x\in \Omega}g(x,\xi)}\over {\xi}}<
\lambda_{1}\ .$$
Assume also that, for each $x\in \Omega$, $g(x,0)=0$ and the function $\xi\to
{{g(x,\xi)}\over {\xi}}$ is non-increasing in $]0,+\infty[$.\par
Let $f:\overline {\Omega}\times {\bf R}\to {\bf R}$ be the function
defined by
$$f(x,\xi)=\cases {g(x,\xi) & if $(x,\xi)\in \overline {\Omega}\times [0,+\infty[$\cr &
\cr 0 & otherwise\ .\cr}$$
Then, the conclusion of Theorem 12 holds.}\par
\smallskip
PROOF.  For each $\lambda>0$, $(x,\xi)\in \overline {\Omega}\times {\bf R}$, put
$$\alpha(\xi)=-(\xi+|\xi|)\xi$$
and
$$h_{\lambda}(x,\xi)=f(x,\xi)+\lambda \alpha(\xi)\ .$$
Clearly, $h_{\lambda}\in {\cal A}$.
Since $h_{\lambda}$ is locally H\"older continuous in $\overline {\Omega}\times {\bf R}$,
the critical points of the functional $\Psi_{h_{\lambda}}$ are continuous in $\overline {\Omega}$.
Thus, since $h_\lambda$ is zero in $\overline {\Omega}\times ]-\infty,0]$, they 
are non-negative in $\Omega$.  Now, observe that, for each $x\in \Omega$,
the function $\xi\to {{h_{\lambda}(x,\xi)}\over {\xi}}$ is (strictly) decreasing in 
$]0,+\infty[$. So, by Theorem 1 of [7], the functional $\Psi_{h_{\lambda}}$ has
at most one non-zero crtitical point in $X$. Consequently, it has at most two
critical points (note that $0$ is one of them). Moreover,
 $\Psi_{h_{\lambda}}$ satisfies the Palais-Smale condition, as it is coercive
(as well as $\Psi_f$).
Thus, since $\Psi_{h_{\lambda}}=\Psi_{f}-\lambda J_{\alpha}$,
all the assumptions
of Theorem 9 are satisfied, and the conclusion follows from it.
\hfill $\bigtriangleup$\par
\medskip
The final result is as follows\par
\medskip
THEOREM 14. - {\it  Let $h=1$, $\Omega=]0,1[$, and let
 $g\in C^{2}([0,+\infty[)$ be a  convex and
non-negative function,
with  $g(0)=0$, such that that
$\sup_{\xi>0}{{g(\xi)}\over
{\xi}}
<\pi^2\ .$
Let $f: {\bf R}\to {\bf R}$ be the function
defined by
$$f(\xi)=\cases {g(\xi) & if $\xi\in  [0,+\infty[$\cr &
\cr 0 & if $\xi<0$\ .\cr}$$
Then, the sublevel sets of $\Psi_f$ are weakly connected.}\par
\smallskip
PROOF. Note that,  in the present case, one has
$\lambda_{1}=\pi^{2}$.  Let $0<\lambda< \pi^{2}-
\sup_{\xi>0}{{g(\xi)}\over {\xi}}$. Define $\alpha:{\bf R}\to
{\bf R}$ by
$$\alpha(\xi)=\cases {\xi-\log(\xi+1) & if $\xi\geq 0$\cr & \cr
0 & if $\xi<0$\ .\cr}$$
Then, the functional $\Psi_{f}-\lambda J_{\alpha}$ satisfies
the Palais-Smale condition (it is coercive as well as $\Psi_f$) and
its crtical points are non-negative.
Clearly, $f+\lambda \alpha\in C^{2}([0,+\infty[)$,
$f(0)+\lambda\alpha(0)=f'(0)+\lambda\alpha'(0)=0$ and
$f''(\xi)+\lambda\alpha''(\xi)>0$ for all $\xi>0$. Then,
by Example 2 of [19], the functional $\Psi_{f}-\lambda J_{\alpha}$ has at most
two critical points in $X$. The conclusion now follows from
Theorem 9.\hfill $\bigtriangleup$\par
\medskip
We conclude with the following problem.\par
\medskip
PROBLEM 2. -  Is there a function $g$ satisfying the hypotheses of
Theorem 14 for which the functional $\Psi_{f}$ has a non-absolute local
minimum ?\par
\vfill\eject
\centerline {{\bf References}}\par
\bigskip
\bigskip
\medskip
\noindent
[1]\hskip 5pt G. ANELLO and G. CORDARO, {\it Existence of solutions
of the Neumann problem involving the p-Laplacian via a variational
principle of Ricceri}, Arch. Math. (Basel), {\bf 79} (2002),
274-287.\par
\smallskip
\noindent
[2]\hskip 5pt G. ANELLO and G. CORDARO, {\it An existence theorem for the
Neumann problem involving the $p$-Laplacian}, J. Convex Anal., {\bf 10}
(2003), 185-198.\par
\smallskip
\noindent
[3]\hskip 5pt G. ANELLO and G. CORDARO, {\it Positive infinitely many and
arbitrarily small solutions for the Dirichlet problem involving the
$p$-Laplacian}, Proc. Royal Soc. Edinburgh Sect. A, {\bf 132} (2002),
511-519.\par
\smallskip
\noindent
[4]\hskip 5pt G. ANELLO, {\it A multiplicity theorem for critical points
of functionals on reflexive Banach spaces}, Arch. Math. (Basel),
to appear.\par
\smallskip
\noindent
[5]\hskip 5pt G. ANELLO, {\it Existence of infinitely many weak solutions
for a Neumann problem}, Nonlinear Anal., to appear.\par
\smallskip
\noindent
[6]\hskip 5pt D. AVERNA and G. BONANNO, {\it A three critical points theorem
and its applications to the ordinary Dirichlet problem}, Topol. Methods
Nonlinear Anal., {\bf 22} (2003), 93-104.\par
\smallskip
\noindent
[7]\hskip 5pt H. BREZIS and L. OSWALD, {\it Remarks on sublinear
elliptic equations}, Nonlinear Anal., {\bf 10} (1986), 55-64.\par
\smallskip
\noindent
[8]\hskip 5pt P. CANDITO, {\it Infinitely many solutions to the Neumann
problem for elliptic equations involving the $p$-Laplacian and with
discontinuous nonlinearities}, Proc. Edinb. Math. Soc., {\bf 45}
(2002), 397-409.\par
\smallskip
\noindent
[9]\hskip 5pt G. CORDARO, {\it Three periodic solutions to an eigenvalue problem
for a class of second order Hamiltonian systems}, Abstr. Appl. Anal.,
{\bf 18} (2003), 1037-1045.\par
\smallskip
\noindent
[10]\hskip 5pt G. CORDARO, {\it Existence and location of periodic solutions
to convex and non coercive Hamiltonian systems}, Discrete Contin. Dyn.
Syst, to appear.\par
\smallskip
\noindent
[11]\hskip 5pt F. FARACI, {\it Bifurcation theorems for
Hammerstein nonlinear integral equations}, Glasg. Math. J., {\bf 44}
(2002), 471-481.\par
\smallskip
\noindent
[12]\hskip 5pt F. FARACI, {\it Multiplicity results for a Neumann problem
involving the $p$-Laplacian}, J. Math. Anal. Appl., {\bf 277} (2003),
180-189.\par
\smallskip
\noindent
[13]\hskip 5pt F. FARACI, {\it Three periodic solutions for a second order
nonautonomous system}, J. Nonlinear Convex Anal., {\bf 3} (2002), 393-399.\par
\smallskip
\noindent
[14]\hskip 5pt F. FARACI and R. LIVREA, {\it Infinitely many periodic
solutions for a second order nonautonomous system}, Nonlinear Anal.,
{\bf 54} (2003), 417-429.\par
\smallskip
\noindent
[15]\hskip 5pt F. FARACI and R. LIVREA, {\it Bifurcation theorems for
nonlinear problems with lack of compactness}, Ann. Polon. Math.,
{\bf 82} (2003), 77-85.\par
\smallskip
\noindent
[16]\hskip 5pt F. FARACI and V. MOROZ, {\it Solutions of Hammerstein
integral equations via a variational principle}, J. Integral Equations
Appl., {\bf 15} (2003), 385-402.\par
\smallskip
\noindent
[17]\hskip 5pt N. GHOUSSOUB and D. PREISS, {\it A general mountain pass principle
for locating and classifying critical points}, Ann. Inst. H. Poincar\'e
Anal. Non Lin\'eaire, {\bf 6} (1989), 321-330.\par
\smallskip
\noindent
[18]\hskip 5pt A. IANNIZZOTTO, {\it A sharp existence and localization
theorem for a Neumann problem}, Arch. Math. (Basel), to appear.\par
\smallskip
\noindent
[19]\hskip 5pt P. KORMAN and Y. LI, {\it Generalized averages for
solutions of two-point Dirichlet problems}, J. Math. Anal. Appl.,
{\bf 239} (1999), 478-484.\par
\smallskip
\noindent
[20]\hskip 5pt S. A. MARANO and D. MOTREANU, {\it Infinitely many critical
points of non-differentiable functions and applications to a Neumann type
problem involving the $p$-Laplacian}, J. Differential Equations, {\bf 182}
(2002), 108-120.\par
\smallskip
\noindent
[21]\hskip 5pt O. NASELLI, {\it A class of functionals on a Banach spaces for
which strong and weak local minima do coincide}, Optimization, {\bf 50}
(2001), 407-411.\par
\smallskip
\noindent
[22]\hskip 5pt B. RICCERI, {\it Some topological mini-max theorems via
an alternative principle for multifunctions}, Arch. Math. (Basel),
{\bf 60} (1993), 367-377.\par
\smallskip
\noindent
[23]\hskip 5pt B. RICCERI, {\it On a topological minimax theorem and
its applications}, in ``Minimax theory and applications'', B. Ricceri
and S. Simons eds., 191-216, Kluwer Academic Publishers, 1998.\par
\smallskip
\noindent
[24]\hskip 5pt B. RICCERI, {\it A general variational principle and
some of its applications}, J. Comput. Appl. Math., {\bf 113}
(2000), 401-410.\par
\smallskip
\noindent
[25]\hskip 5pt B. RICCERI, {\it Infinitely many solutions of the Neumann
problem for elliptic equations involving the p-Laplacian},
Bull. London Math. Soc., {\bf 33} (2001), 331-340.\par
\smallskip
\noindent
[26]\hskip 5pt B. RICCERI, {\it On a classical existence
theorem for nonlinear elliptic equations}, in ``Experimental,
constructive and nonlinear analysis'', M. Th\'era ed.,
275-278, CMS Conf. Proc. {\bf 27}, Canad. Math. Soc., 2000.\par
\smallskip
\noindent
[27]\hskip 5pt B. RICCERI, {\it A bifurcation theory for some nonlinear
elliptic equations}, Colloq. Math., {\bf 95} (2003), 139-151.\par
\smallskip
\noindent
[28]\hskip 5pt J. SAINT RAYMOND, {\it On the multiplicity of the solutions
of the equation $-\Delta u=\lambda f(u)$}, J. Differential Equations,
{\bf 180} (2002), 65-88.\par
\smallskip
\noindent
[29]\hskip 5pt E. ZEIDLER, {\it Nonlinear functional analysis and
its applications}, vol. III, Springer-Verlag, 1985.\par

\bigskip
\bigskip
\bigskip
\bigskip
Department of Mathematics\par
University of Catania\par
Viale A. Doria 6\par
95125 Catania\par
Italy\par
{\it e-mail}: ricceri@dmi.unict.it
\bye